\documentclass{article}
\usepackage{amsmath}
\usepackage{amssymb}
\usepackage[active]{srcltx} 
\usepackage{amssymb,amsmath,amsthm}
\usepackage{bm}
\usepackage{multicol}
\usepackage{tikz}
\usetikzlibrary{calc}
\RequirePackage{hyperref}
\usepackage{libertinus}\providecommand{\coloneqq}{:=}
\theoremstyle{plain}
\newtheorem{theo}{Theorem}[section]

\newtheorem{lemm}[theo]{Lemma}
\newtheorem*{lemm*}{Lemma}

\newtheorem{namet}[theo]{\myThmName}
\newtheorem*{namet*}{\myThmName}

\newenvironment*{nthm*}[1][\kern-.35em]{\edef\myThmName{#1}\begin{namet*}}{\end{namet*}}

\theoremstyle{definition}

\newtheorem*{defi*}{Definition}

\newtheorem{defs}[theo]{Definitions}

\newtheorem{exas}[theo]{Examples}

\newtheorem*{rema*}{Remark}

\newtheorem{named}[theo]{\myThmName}
\newtheorem*{named*}{\myThmName}

\newcounter{rememberEnumi}

\let\setminus\smallsetminus
\newcommand{\FF}{\ensuremath{\mathbb{F}}}

\let\setminus\smallsetminus
\makeatletter
\newcommand{\gal}[1]{ {\mathchoice
    {\vbox
    {\m@th \ialign {##\crcr \noalign {\kern 1\p@ }\kern 1\p@ \hrulefill \crcr
        \noalign {\kern 1\p@ \nointerlineskip }%
        $\hfil \displaystyle {#1}\hfil $\crcr }}}
    {\vbox
    {\m@th \ialign {##\crcr \noalign {\kern 1\p@ }\kern 1\p@ \hrulefill \crcr
        \noalign {\kern 1\p@ \nointerlineskip }%
        $\hfil \textstyle {#1} $\crcr }}}
    {\vbox
    {\m@th \ialign {##\crcr \noalign {\kern 1\p@ }\kern 1\p@ \hrulefill \crcr
        \noalign {\kern 1\p@ \nointerlineskip }%
        $\hfil \scriptstyle {#1}\hfil $\crcr }}}
    {\vbox
    {\m@th \ialign {##\crcr \noalign {\kern 1\p@ }\kern 1\p@ \hrulefill \crcr
        \noalign {\kern 1\p@ \nointerlineskip }%
        $\hfil \scriptscriptstyle {#1}\hfil $\crcr }}}%
    }}
\makeatother
\newcommand{\set}[2]{\left\{{#1}\left|\vphantom{#1#2\strut}\right.\, 
                    {#2}\right\}}
\newcommand{\Aut}[2][]{\operatorname{Aut}_{#1}(#2)}

\newcommand{\gU}[2]{\operatorname{\Gamma U}(#1,#2)}
\newcommand{\PSU}[2]{\operatorname{PSU}(#1,#2)}
\newcommand{\PgU}[2]{\operatorname{P\Gamma U}(#1,#2)}
\newcommand{\cB}{\mathcal{B}}
\newcommand{\cH}{\mathcal{H}}
\newcommand{\cU}{\mathcal{U}}
\newcommand{\MSC}[1]{\href{https://mathscinet.ams.org/mathscinet/freetools/msc-search?text=#1}{#1}}
  \title{Unitals without O'Nan configurations are classical if they
    admit all translations} %
  \author{Markus J.~Stroppel}%
\begin{document}
%%%%%%%%%%%%%%%%%%%%%%%%%%%%%%%%%%%%%%%%%%%%%%%%%%%%%%%%%%%%%%%%%%%%%%
\maketitle
%%%%%%%%%%%%%%%%%%%%%%%%%%%%%%%%%%%%%%%%%%%%%%%%%%%%%%%%%%%%%%%%%%%%%%
%%%%%%%%%%%%%%%%%%%%%%%%%%%%%%%%%%%%%%%%%%%%%%%%%%%%%%%%%%%%%%%%%%%%%%
  \begin{abstract}\noindent %
    We prove the statement in the title: if a (finite) unital admits
    all translations and contains no O'Nan configurations then the
    unital is classical, i.e., isomorphic to the Hermitian unital of
    the same order.


    MSC 2020 classification:%
    \MSC{51E05} %% General block designs in finite geometry
    \MSC{51A45} %% Incidence structures embeddable into projective geometries
    \MSC{05B05} %% Combinatorial aspects of block designs
  \end{abstract}

\section*{Introduction}

We attempt a brief overview of the background for the present note.
See Section~\ref{sec:unitals} for precise definitions. 

In~\cite{MR690826}, Wilbrink has characterized the finite classical
unitals (i.e., those defined by a Hermitian form,
see~\ref{exas:classical} below) among all finite unitals by three
elementary conditions (see~\ref{wilbrink} below). We show that these
conditions are satisfied if the unital in question admits all
translations, and contains no O'Nan configurations.  Then Wilbrink's
characterization gives the result announced in the title.

Unitals admitting all translations have been studied
in~\cite{MR3090721}. Using the classification of finite simple groups,
it is proved in that paper that such unitals are classical. %
Our present treatment is much more elementary, but imposes the fairly
strong hypothesis that no O'Nan configurations appear in the unital. %
Using translations (assumed to exist by our second strong hypothesis)
we then check Wilbrink's conditions (see~\ref{wilbrink}), and obtain
the result.

It has been conjectured (see~\cite[p.\,102]{MR587626},
\cite[p.\,87]{MR2440325}) that finite classical unitals are
characterized by the absence of O’Nan configurations. To the author's
knowledge, this conjecture has not been proved yet. In many of the
known non-classical unitals, O'Nan configurations have been found.

\section{Unitals, translations}
\label{sec:unitals}

\enlargethispage{4mm}%
\begin{defs}
  Throughout the present paper, let $\cU = (P,\cB,\in)$ be a finite
  \emph{unital} of order~$q$; i.e., the set~$P$ of points has size
  $q^3+1$ with $q>1$, the set~$\cB$ of blocks consists of subsets of
  size~$q+1$ in~$P$, and each $2$-element subset $\{x,y\}\subset P$ is
  contained in exactly one block (called the \emph{joining} block, and
  denoted by $x\vee y$). %
  In other words, $\cU$ is a $2$-$(q^3+1,q+1,1)$-design. 

\goodbreak%
  An \emph{O'Nan configuration} is a set of~$4$ mutually intersecting
  blocks together with~$6$ points, such that each one of the~$6$
  points lies on exactly two of these blocks. See
  Figure~\ref{fig:ONan}. (In axiomatic projective geometry, that
  configuration is called Veblen-Young configuration.)

  \begin{figure}[h!]
    \centering
    \begin{tikzpicture}[%
      yscale=.7, every node/.append style={circle, fill=black, %
        inner sep=.5pt, %
        minimum size=3pt}]%
      \node[coordinate] (Y) at (0,0) {} ;%
      \node[coordinate] (X) at (0,2) {} ;%
      \node[coordinate] (W) at (5,.5) {} ;%
      \node[coordinate] (Z) at ($(Y)!.7!(W)$) {} ;%
      \node[coordinate] (Zt) at ($(X)!.7!(Z)$) {} ;%
      \node[coordinate] (Zx) at ($(X)!1.2!(Z)$) {} ;%
      \node[coordinate] (Wx) at ($(X)!1.2!(W)$) {} ;%
      \node[coordinate] (Wt) at (intersection of X--W and Y--Zt) {} ;%
      \node[coordinate] (Dx) at ($(Y)!1.2!(W)$) {} ;%
      \node[coordinate] (Dtx) at ($(Y)!1.4!(Wt)$) {} ;%

      \draw (X) -- (Zx) ;%
      \draw (X) -- (Wx) ;%
      \draw (Y) -- (Dx) ;%
      \draw (Y) -- (Dtx) ;%

      \node at (X) {} ;%
      \node at (Y) {} ;%
      \node at (Z) {} ;%
      \node at (W) {} ;%
      \node at (Zt) {} ;%
      \node at (Wt) {} ;%
    \end{tikzpicture}
    \caption{An O'Nan configuration.}
    \label{fig:ONan}
  \end{figure}
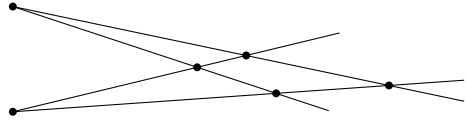

  A \emph{translation} (with center~$c$) of the
  unital~$\cU$ is an automorphism of~$\cU$ that fixes a
  point~$c$ and every block through~$c$.
\end{defs}

\begin{lemm}[\protect{\cite[Thm.\,1.3]{MR3090721}}]
  For each point~$c$ of\/~$\cU$, the set\/ $T_{[c]}$ of all
  translations with center~$c$ forms a subgroup of\/ $\Aut\cU$ that
  acts semi-regularly on $P\setminus\{c\}$. %
\end{lemm}

In particular, the group $T_{[c]}$ induces a semi-regular group on
$B\setminus\{c\}$, for each block $B$ containing~$c$.  We say that
$\cU$ \emph{admits all translations with center~$c$} if $T_{[c]}$ is
transitive on $B\setminus\{c\}$. (This is equivalent to
$|T_{[c]}| = q$.)  If $\cU$ admits all translations with center~$c$,
for each $c\in P$, we say that~$\cU$ \emph{admits all translations}.
In~\cite{MR3090721}, it has been proved that each unital admitting all
translations is isomorphic to the classical unital of the same order
(see~\ref{exas:classical} below). %
That proof uses the classification of finite simple groups. %

There do exist examples of unitals admitting all translations for many
(but not all) centers (see~\cite{MR3533345}, \cite{MR4581128}), and
unitals admitting all translations for just one single center
(see~\cite[Sec.\,5]{MR3533345}).  At least some of those unitals do
contain O'Nan configurations (e.g., see~\cite[6.7,
6.10]{Moehler2020}). %

\begin{exas}\label{exas:classical}
  Let $E|F$ be a separable quadratic field extension, and
  let~$\sigma\colon s\mapsto\gal{s}$ denote the generator of the
  Galois group. There is an (essentially unique) isotropic
  non-degenerate $\sigma$-Hermitian form~$h$ on~$E^3$.  The classical
  unital $\cH_{E|F}$ (cp. \cite[p.\,104]{MR0233275}, %
  \cite[2.1, 2.2, see also p.\,29]{MR2440325}) has the set $P_{E|F}$
  of all one-dimensional isotropic subspaces of~$E^3$ as point set,
  the blocks are the intersections of $P_{E|F}$ with lines that meet
  that set in more than one point. %
  Up to a choice of basis, the Hermitian form~$h$ maps
  $(x,y)\in E^3\times E^3$ with $x=(x_0,x_1,x_2)$, $y=(y_0,y_1,y_2)$
  to $x_0\gal{y_2}-x_1\gal{y_1}+x_2\gal{y_0}$. %

  We note that the classical unital $\cU_{E|F}$ does not contain any
  O'Nan configurations. %
  (See~\cite[Proposition, p.\,507]{MR0295934} for the finite case,
  and~\cite[2.2]{MR2795696}). 

  The automorphism group $\Aut{\cH_{E|F}}$ of the classical
  unital~$\cH_{E|F}$ is the group $\PgU{E^3}{h}$ induced by the group
  $\gU{E^3}h$ of all semi-similitudes of the form~$h$.  %
  (See~\cite{MR0295934} for the finite case, and~\cite{Tits} for the
  general case; cp.~\cite[6.1, 5.6]{MR2241352}.) %
  If the form is given as above, the subgroup
  \[
   \Xi \coloneqq \set{\left(
        \begin{matrix}
          1 & x & z \\
          0 & 1 & \gal{x} \\
          0 & 0 & 1
        \end{matrix}\right)}
    {x,z\in E, z+\gal{z}=x\gal{x}}
  \]
  is contained in the stabilizer of the point $E(0,0,1) \in P_{E|F}$,
  and
  \[
    T_{[E(0,0,1)]} = \set{\left(
        \begin{matrix}
          1 & 0 & z \\
          0 & 1 & 0 \\
          0 & 0 & 1
        \end{matrix}\right)}
    {z\in E, z+\gal{z}=0} 
  \]
  is the group of translations with center~$E(0,0,1)$.  The
  group~$\Xi$, together with any element of $\PSU{E^3}h$ that
  moves~$E(0,0,1)$, shows that $\PSU{E^3}h \le \Aut{\cH_{E|F}}$ acts
  doubly transitively on~$P_{E|F}$.

  In particular, if $E$ is finite of square order~$q^2$ then there is
  a unique subfield~$F$ of order~$q$. The involution~$\sigma$ is the
  appropriate power of the Frobenius endomorphism, it maps $s\in E$
  to~$\gal{s} = s^q$.  Then $\cH_{E|F} \cong \cH_{\FF_{q^2}|\FF_q}$ is
  a unital of order~$q$; the groups considered above show that this
  unital admits all translations.
\end{exas}

\section{Wilbrink's conditions}

\begin{defs}
  Consider $B\in\cB$ and $x\in P\setminus B$. A block $B'$ with
  $x\notin B'$ is called \emph{$x$-parallel to~$B$} if $B'$ meets each
  block joining~$x$ with a point on~$B$.

  If\/~$\cU$ does not contain any O'Nan configurations then there
  exists at most one block through a given point~$y\ne x$ that is
  $x$-parallel to a given block~$B$. %
  If $\tau$ is a translation with center~$x$ then the image~$B^\tau$
  of~$B$ under~$\tau$ is $x$-parallel to~$B$.
\end{defs}
    
\begin{theo}[\protect{\cite{MR690826}}]\label{wilbrink}%
  Let\/ $\cU = (P,\cB,\in)$ be a unital of order~$q$ satisfying the
  following conditions.
  \begin{itemize}
  \item[\upshape(I)]%
    There are no O'Nan configurations in~$\cU$. 
  \item[\upshape(II)]%
    Consider $L\in\cB$ and\/ $x,y\in P$ such that\/ $x\notin L$ and\/
    $(x\vee y)\cap B \ne\emptyset$. Then there exists an $x$-parallel
    block\/~$L$ through~$y$.
  \item [\upshape(III)]%
    Consider three blocks $M_0,M_1,M_2$ through a common point~$x$, and points
    $y_i,z_i\in M_i\setminus\{x\}$ for $i\in\{0,1,2\}$. %
    If $z_0\vee z_i$ is $x$-parallel to $y_0\vee y_i$ for
    $i\in\{1,2\}$ then $z_1\vee z_2$ is $x$-parallel to $y_1\vee y_2$.
  \end{itemize}
  Then $\cU$ is isomorphic to the classical unital of order~$q$. %
\goodbreak%
\end{theo}

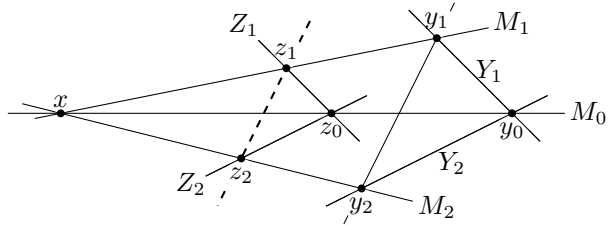
\begin{figure}[h]
  \centering
  \begin{tikzpicture}[%
    scale=1, %
    every node/.append style={circle, fill=black, %
      inner sep=.5pt, %
      minimum size=3pt}, %
    every path/.append style={shorten <=-15pt, shorten >=-15pt} %
    ]%
    \node[coordinate] (x) at (0,0) {} ;%
    \node[coordinate] (y0) at (6,0) {} ;%
    \node[coordinate] (y1) at (5,1) {} ;%
    \node[coordinate] (y2) at (4,-1) {} ;%
    \node[coordinate] (z0) at ($(y0)!.4!(x)$) {} ;%
    \node[coordinate] (z1) at ($(y1)!.4!(x)$) {} ;%
    \node[coordinate] (z2) at ($(y2)!.4!(x)$) {} ;%

    \draw[shorten <=-20pt, shorten >=-20pt] (x) -- (y0) ;%
    \draw[shorten <=-10pt, shorten >=-20pt] (x) -- (y1) ;%
    \draw[shorten <=-15pt, shorten >=-20pt] (x) -- (y2) ;%

    \draw (y0) -- (y1) ;%
    \draw (y0) -- (y2) ;%
    \draw (y1) -- (y2) ;%
    \draw (z0) -- (z1) ;%
    \draw (z0) -- (z2) ;%
    \draw[dashed, thick, shorten <=-20pt, shorten >=-20pt] (z1) -- (z2) ;%

    \node[fill=white, draw=none, anchor=south] at (x)   {\(x\) } ;%
    \node[fill=white, draw=none, anchor=north] at (y0)  {\(y_0\)} ;%
    \node[fill=white, draw=none, anchor=south] at (y1)  {\(y_1\)} ;%
    \node[fill=white, draw=none, anchor=north] at (y2)  {\(y_2\)} ;%
    \node[fill=white, draw=none, anchor=north] at (z0)  {\(z_0\)} ;%
    \node[fill=white, draw=none, anchor=south] at (z1)  {\(z_1\)} ;%
    \node[fill=white, draw=none=north, anchor=north] at (z2)  {\(z_2\)} ;%

    \node[fill=none, draw=none, anchor=south] at ($(y1)!.7!(y0)$) {$Y_1$} ;%
    \node[fill=none, draw=none, anchor=north] at ($(y2)!.6!(y0)$) {$Y_2$} ;%
    \node[fill=none, draw=none] at ($(z1)!-.95!(z0)$) {$Z_1$} ;%
    \node[fill=none, draw=none] at ($(z2)!-.56!(z0)$) {$Z_2$} ;%
    \node[fill=none, draw=none] at ($(x)!1.17!(y0)$) {$M_0$} ;%
    \node[fill=none, draw=none] at ($(x)!1.20!(y1)$) {$M_1$} ;%
    \node[fill=none, draw=none] at ($(x)!1.25!(y2)$) {$M_2$} ;%

    \node at (x)   {} ;%
    \node at (y0)  {} ;%
    \node at (y1)  {} ;%
    \node at (y2)  {} ;%
    \node at (z0)  {} ;%
    \node at (z1)  {} ;%
    \node at (z2)  {} ;%

    \draw (y0) -- (y1) ;%
    \draw (y0) -- (y2) ;%
    \draw (z0) -- (z1) ;%
    \draw (z0) -- (z2) ;%

  \end{tikzpicture}
  \caption{Wilbrink's condition~(III)}
  \label{fig:wilbrinkIII}
\end{figure}

\begin{theo}\label{theTheo}
  Let\/ $\cU = (P,\cB,\in)$ be a unital of order~$q$ admitting all
  translations, and with no O'Nan configurations. Then $\cU$ is
  isomorphic to the classical unital of order~$q$.
\end{theo}
\begin{proof}
  It remains to show that~$\cU$ satisfies Wilbrink's conditions~(II)
  and~(III), see~\ref{wilbrink}.

  Consider $B\in\cB$ and $x\in P\setminus B$. If $x\vee y$ meets~$B$
  in a point~$w$, say, then our assumptions secure the existence of a
  translation $\tau\in T_{[x]}$ mapping $w$ to~$y$. Clearly the
  image~$B'$ of~$B$ under~$\tau$ meets each block through~$x$ that
  also meets~$B$. So $B'$ is $x$-parallel to~$B$, and condition~(II)
  is verified.

  Now consider three blocks $M_0,M_1,M_2$ through~$x$, and points
  $y_i,z_i\in M_i\setminus\{x\}$ for $i\in\{0,1,2\}$ %
  (see Figure~\ref{fig:wilbrinkIII}). %
  Let $\tau$ be the translation with center~$x$ that maps~$y_0$
  to~$z_0$. %
  For $i\in\{1,2\}$, the image~$Z_i$ of $Y_i \coloneqq y_0\vee y_i$
  under~$\tau$ contains~$z_0$, and meets each block through~$x$ that
  meets $Y_i$. So $Z_i$ is $x$-parallel to~$Y_i$. By the absence of
  O'Nan configurations, we know that the $x$-parallel block to~$Y_i$
  through~$z_0$ is unique, and infer $Z_i = z_0\vee z_i$. Thus~$z_i$
  is the image of~$y_i$ under~$\tau$, and $z_1\vee z_2$ is the image
  of $y_1\vee y_2$ under~$\tau$. This yields that $z_1\vee z_2$ is
  $x$-parallel to $y_1\vee y_2$, and condition~(III) is verified.  
\end{proof}

An application of the present theorem will be given
in~\cite{StroppelInPrep}.

\bigskip\noindent%
M.J. Stroppel\\[1ex]
LExMath \\
Fachbereich Mathematik\\
Universit\"at Stuttgart \\
70550 Stuttgart \\
Germany %
\\[1ex]
    {stroppel@mathematik.uni-stuttgart.de}%

\end{document}